\documentclass[12pt]{amsart}

\usepackage{amsfonts,amsmath,amssymb,amsthm}

\newcommand{\Cc}{\mathbb{C}}

\newcommand{\Zz}{\mathbb{Z}}
\newcommand{\defi}[1]{\emph{#1}}

\renewcommand{\epsilon}{\varepsilon}
\renewcommand{\le}{\leqslant}
\renewcommand{\ge}{\geqslant}

\newcommand {\grad}{\mathop{\mathrm{grad}}\nolimits}
\newcommand {\ord}{\mathop{\mathrm{ord}}\nolimits}
\newcommand {\mult}{\mathop{\mathrm{mult}}\nolimits}
\newcommand {\xubar}{\underline{x}}
\newcommand {\charac}{\mathop{\mathrm{char}}\nolimits}

{\theoremstyle{plain}
\newtheorem{theorem}{Theorem}[section]    
\newtheorem{lemma}[theorem]{Lemma}       
\newtheorem{proposition}[theorem]{Proposition}      
\newtheorem{corollary}[theorem]{Corollary}      
}
{\theoremstyle{remark}
\newtheorem{definition}[theorem]{Definition}      
\newtheorem*{remark*}{Remark}  
\newtheorem{remark}[theorem]{Remark}   

}

\begin{document}
\title{Reducibility of rational functions in several variables}

\author{Arnaud Bodin}
\address{Laboratoire Paul Painlev\'e, Math\'ematiques,
Universit\'e de Lille 1,  59655 Villeneuve d'Ascq, France.}
\email{Arnaud.Bodin@math.univ-lille1.fr}

\date{\today}

\begin{abstract}
  We prove a analogous of Stein theorem for rational functions in
  several variables: we bound the number of reducible fibers by a
  formula depending on the degree of the fraction.
\end{abstract}

\maketitle

\section{Introduction}

Let $K$ be an algebraically closed field.  Let $f = \frac pq \in
K(\xubar)$, with $\xubar = (x_1,\ldots,x_n)$, $n\ge 2$ and
$\gcd(p,q)=1$, the \defi{degree} of $f$ is $\deg f = \max \{ \deg p,
\deg q\}$.  We associate to a fraction $f = \frac pq$ the pencil
$p-\lambda q$, $\lambda \in \hat K$ (where we denote $\hat K = K \cup
\{\infty\}$ and by convention if $\lambda = \infty$ then $p-\lambda q=
q$).

For each $\lambda \in \hat K$ write the decomposition into irreducible
factors:
$$p-\lambda q = \prod_{i=1}^{n_\lambda} F_{i}^{r_i}.$$
The
\defi{spectrum} of $f$ is $\sigma(f)=\{ \lambda \in \hat K \mid
n_\lambda >1\}$, and the \defi{order of reducibility} is $\rho(f) =
\sum_{\lambda \in \hat K} (n_\lambda-1).$

A fraction $f$ is \defi{composite} if it is the composition of a
univariate rational fraction of degree more than $1$ with another
rational function.

\begin{theorem}
  Let $K$ be an algebraically closed field of characteristic $0$. 
  Let $f \in K(\xubar)$ be non-composite then
  $$\rho(f) < (\deg f)^2+\deg f.$$
\end{theorem}

A theorem of Bertini and Krull implies that if $f$ is non-composite
then $\sigma(f)$ is finite and we should notice that $\# \sigma(f) \le
\rho(f)$.  Later on, for an algebraically closed field of
characteristic zero and for a polynomial $f\in K[x,y]$, Stein
\cite{St} proved the formula $\rho(f) < \deg f$.  This formula has
been generalized in several directions, see \cite{Na} for references.
For a rational function $f\in \Cc(x,y)$ a consequence of the work of
Ruppert \cite{Ru} on pencil of curves, is that $\# \sigma(f) < (\deg f)^2$.  For $K$
algebraically closed (of any characteristic) and $f\in K(x,y)$
Lorenzini \cite{Lo} proved under geometric hypotheses 
on the pencil $(p-\lambda q)$ that $\rho(f) < (\deg f)^2$. 
This has been generalized 
by Vistoli \cite{Vi} for a pencil in several variables
for an algebraically closed field of characteristic $0$.

Let us give an example extracted from \cite{Lo}. Let $f(x,y) =
\frac{x^3+y^3+(1+x+y)^3}{xy(1+x+y)}$, then $\deg(f)=3$ and $\sigma(f)=
\{1,j,j^2,\infty\}$ (where $\{1,j,j^2\}$ are the third roots of
unity).  For $\lambda \in \sigma(f)$, $(f=\lambda)$ is composed of
three lines hence $\rho(f) = 8 = (\deg f)^2-1$. Then Lorenzini's bound
is optimal in two variables.

The motivation of this work is that we develop the analogous theory of
Stein for rational function: composite fractions, kernels of Jacobian
derivatives, groups of divisors,...   The method for the two
variables case is inspired from the work of Stein \cite{St} and the
presentation of that work by Najib \cite{Na}.  For completeness even
the proofs similar to the ones of Stein have been included.  Another
motivation is that with a bit more effort we get the case of several
variables by following the ideas of \cite{Na} (see the articles
\cite{Na1}, \cite{Na2}).

In \S \ref{sec:composite} we prove that a fraction is non-composite if
and only its spectrum is finite.  Then in \S \ref{sec:kernel} we
introduce a theory of Jacobian derivation and compute the kernel.
Next in \S \ref{sec:order} we prove that for a non-composite fraction
in two variables $\rho(f) < (\deg f)^2+\deg f$.  Finally in \S
\ref{sec:several} we extend this formula to several variables and we
end by stating a result for fields of any characteristic.

\emph{Acknowledgements:} I wish to thank Pierre D\`ebes and Salah
Najib for discussions and encouragements.

\section{Composite rational functions}
\label{sec:composite}

Let $K$ be an algebraically closed field.  Let $\xubar =
(x_1,\ldots,x_n)$, $n\ge 2$.

\begin{definition}
  A rational function $f \in K(\xubar)$ is \defi{composite} if there
  exist $g\in K(\xubar)$ and $r\in K(t)$ with $\deg r \ge 2$ such that
  $$f = r\circ g.$$
\end{definition}

\begin{theorem}
\label{th:composite}
Let $f = \frac pq \in K(\xubar)$. The following assertions are
equivalent:
\begin{enumerate}
\item \label{it:i} $f$ is composite;
\item \label{it:ii} $p-\lambda q$ is reducible in $K[\xubar]$ for all
  $\lambda \in \hat K$ such that $\deg p-\lambda q= \deg f$;
\item \label{it:iii} $p-\lambda q$ is reducible in $K[\xubar]$ for
  infinitely many $\lambda \in \hat K$.
\end{enumerate}
\end{theorem}

Before proving this result we give two corollaries.

\begin{corollary}
\label{cor:composite}
$f$ is non-composite if and only if its spectrum $\sigma(f)$ is
finite.
\end{corollary}

One aim of this paper is to give a bound for $\sigma(f)$.  The hard
implication of this theorem $(\ref{it:iii}) \Rightarrow (\ref{it:i})$
is in fact a reformulation of a theorem of Bertini and Krull.

We also give a nice application pointed out to us by P.~D\'ebes:
\begin{corollary}
Let $p\in K[\xubar]$ irreducible. Let $q\in K[\xubar]$ with
$\deg q < \deg p$ and $\gcd(p,q)=1$. Then for all but finitely 
many $\lambda \in K$, 
$p-\lambda q$ is irreducible in $K[\xubar]$.
\end{corollary}

\medskip

\emph{Convention :} When we define a fraction $F = \frac PQ$ we will
assume that $\gcd(P,Q)=1$.

\medskip

We start with the easy part of Theorem \ref{th:composite}:
\begin{proof}
  (\ref{it:ii}) $\Rightarrow$ (\ref{it:iii}) is trivial.  Let us prove
  (\ref{it:i}) $\Rightarrow$ (\ref{it:ii}).  Let $f=\frac pq$ be a
  composite rational function.  There exist $g=\frac{u}{v}\in
  K(\xubar)$ and $r\in K(t)$ with $k=\deg r \ge 2$ such that $ f =
  r\circ g$ . Let us write $r=\frac ab$. Let $\lambda \in \hat K$ such that
  $\deg a- \lambda b = \deg r$ and factorize $a(t)-\lambda b(t)=
  \alpha{(t-t_1)(t-t_2)\cdots(t-t_k)}$, $\alpha\in K^*, t_1,\ldots,t_k
  \in K$.  Then
  $$p-\lambda q = q \cdot \left(f-\lambda\right) = q\cdot
  \left(\frac{a-\lambda b}{b}\right)(g) = \alpha q
  \frac{(g-t_1)\cdots(g-t_k)}{b(g)}.$$
  Then by multiplication by $v^k$
  at the numerator and denominator we get:
  $$(p-\lambda q)\cdot (v^k b(g)) = \alpha q (u-t_1 v)\cdots(u-t_k
  v),$$
  which is a polynomial identity. As $\gcd(a,b)=1$,
  $\gcd(u,v)=1$ and $\gcd(p,q)=1$ then $u-t_1 v,\ldots,u-t_k v$ divide
  $p-\lambda q$.  Hence $p-\lambda q$ is reducible in $K[\xubar]$.
 \end{proof}

Let us reformulate the Bertini-Krull theorem in our context from
\cite[Theorem 37]{Sc}. It will enable us to end the proof of Theorem 
\ref{th:composite}.
\begin{theorem}[Bertini, Krull]
  Let $F(\xubar,\lambda) = p(\xubar)-\lambda q(\xubar) \in
  K[\xubar,\lambda]$ an irreducible polynomial.  Then the following
  conditions are equivalent:
\begin{enumerate}
\item $F(\xubar,\lambda_0) \in K[\xubar]$ is reducible for all
  $\lambda_0 \in K$ such that $\deg_{\xubar} F(\xubar,\lambda_0) =
  \deg_{\xubar} F$.
\item
  \begin{enumerate}
  \item either there exist $\phi,\psi \in K[\xubar]$ with
    $\deg_{\xubar} F > \max\{\deg \phi, \deg \psi\}$, and $a_i\in
    K[\lambda]$, such that
    $$F(\xubar,\lambda) = \sum_{i=0}^n
    a_i(\lambda)\phi(\xubar)^{n-i}\psi(\xubar)^{i} ; $$
  \item or $\charac(K)=\pi > 0$ and $F(\xubar,\lambda) \in
    K[\xubar^\pi,\lambda]$, where $\xubar^\pi =
    (x_1^\pi,\ldots,x_n^\pi)$.
  \end{enumerate}
\end{enumerate}
\end{theorem}

We now end the proof of Theorem \ref{th:composite}:
\begin{proof}
  (\ref{it:iii}) $\Rightarrow$ (\ref{it:i}) Suppose that $p-\lambda_0
  q$ is reducible in $K[\xubar]$ for infinitely many $\lambda_0 \in
  \hat K$; then it is reducible for all $\lambda_0 \in K$ such that
  $\deg_{\xubar} F(\xubar,\lambda_0) = \deg_{\xubar} F$ (see Corollary
  3 of Theorem 32 of \cite{Sc}).  We apply Bertini-Krull theorem:
  
  \emph{Case (a):} $F(\xubar,\lambda) = p(\xubar)-\lambda q(\xubar)$
  can be written:
  $$
  p(\xubar)-\lambda q(\xubar) = \sum_{i=0}^n
  a_i(\lambda)\phi(\xubar)^{n-i}\psi(\xubar)^{i}.$$
  So we may suppose
  that for $i=1,\ldots,n$, $\deg_\lambda a_i =1$, let us write
  $a_i(\lambda)= \alpha_i-\lambda\beta_i$, $\alpha_i,\beta_i\in K$.
  Then
  $$p(\xubar) = \sum_{i=0}^n \alpha_i
  \phi(\xubar)^{n-i}\psi(\xubar)^{i}= \phi^n \sum_{i=0}^n \alpha_i
  \Big(\frac{\psi}{\phi}\Big)^i(\xubar), $$
  and
  $$
  q(\xubar) = \sum_{i=0}^n
  \beta_i\phi(\xubar)^{n-i}\psi(\xubar)^{i}= \phi^n \sum_{i=0}^n
  \beta_i \Big(\frac{\psi}{\phi}\Big)^i(\xubar).$$
  
  If we set $g(\xubar)= \frac{\psi(\xubar)}{\phi(\xubar)} \in
  K[\xubar]$, and $r(t) = \frac{\sum_{i=0}^n \alpha_i
    t^i}{\sum_{i=0}^n \beta_i t^i}$ then $\frac{p}{q}(\xubar)= r \circ
  g$.  Moreover as $\deg_{\xubar} F > \max \{ \deg\phi, \deg \psi\}$
  this implies $n\ge 2$ so that $\deg r \ge 2$.  Then $ \frac pq= f =
  r\circ g$ is a composite rational function
  
  \emph{Case (b):} Let $\pi = \charac(K) >0$ and
  $F(\xubar,\lambda)=p(\xubar)-\lambda q(\xubar) \in
  K[\xubar^\pi,\lambda]$, For $\lambda = 0$ it implies that $p(\xubar)
  = P(\xubar^\pi)$, then there exists $p'\in K[\xubar]$ such that
  $p(\xubar) = \left( p'(\xubar) \right)^\pi$. For $\lambda = -1$ we
  obtain $s'\in K[\xubar]$ such that $p(\xubar)+q(\xubar) = \left(
    s'(\xubar) \right)^\pi$. Then $q(\xubar)= (p(\xubar)+q(\xubar)) -
  p(\xubar)= \left( s'(\xubar) \right)^\pi - \left( p'(\xubar)
  \right)^\pi = \left( s'(\xubar)-p'(\xubar) \right)^\pi$.  Then if we
  set $q'=s'-p'$ we obtain $q(\xubar) = \left( q'(\xubar)
  \right)^\pi$.  Now set $r(t) = t^\pi$ and $g=\frac {p'}{q'}$ we get
  $f = \frac pq= \left(\frac{p'}{q'}\right)^\pi = r\circ g$.
\end{proof}

\section{Kernel of the Jacobian derivation}
\label{sec:kernel}

We now consider the two variables case and $K$ is an uncountable
algebraically closed field of characteristic zero.

\subsection{Jacobian derivation}

Let $f,g \in K(x,y)$, the following formula:
$$D_f(g) = \frac{\partial f}{\partial x}\frac{\partial g}{\partial y}-
\frac{\partial f}{\partial y}\frac{\partial g}{\partial x},$$
defines
a derivation $D_f : K(x,y) \rightarrow K(x,y).$ Notice the $D_f(g)$ is
the determinant of the Jacobian matrix of $(f,g)$.  We denote by $C_f$
the kernel of $D_f$:
$$C_f = \left\lbrace g \in K(x,y) \mid D_f(g)=0 \right\rbrace.$$
Then
$C_f$ is a subfield of $K(x,y)$. We have the inclusion $K(f) \subset
C_f$.  Moreover if $g^k\in C_f$, $k\in \Zz\setminus\{0\}$ then $g\in
C_f$.

\begin{lemma}
\label{lem:kernel}
Let $f=\frac pq$, $g \in K(x,y)$. The following conditions are
equivalent:
\begin{enumerate}
\item \label{it:li} $g \in C_f$;
\item \label{it:lii} $f$ and $g$ are algebraically dependent;
\item \label{it:liii} $g$ is constant on irreducible components of the
  curves $(p-\lambda q=0)$ for all but finitely many $\lambda \in \hat
  K$;
\item \label{it:liv} $g$ is constant on infinitely many irreducible
  components of the curves $(p-\lambda q=0)$, $\lambda \in \hat K$.
\end{enumerate}
\end{lemma}

\begin{corollary}
\label{cor:kernel}
If $g\in C_f$ is not a constant then $C_f = C_g$.
\end{corollary}

\begin {proof}\ 
\begin{itemize}
\item (\ref{it:li}) $\Leftrightarrow$ (\ref{it:lii}).  We follow the
  idea of \cite{Na} instead of \cite{St}.  $f$ and $g$ are
  algebraically dependent if and only $\mathrm{transc}_K K(f,g) =1$.
  And $\mathrm{transc}_K K(f,g) =1$ if and only the rank of the
  Jacobian matrix of $(f,g)$ is less or equal to $1$, which is
  equivalent to $g \in C_f$.
  
\item (\ref{it:lii}) $\Rightarrow$ (\ref{it:liii}).  Let $f$ and $g$
  be algebraically dependent. Then there exists a two variables
  polynomial in $f$ and $g$ that vanishes. Let us write
  $$\sum_{i=0}^{n} R_i(f)g^i=0$$
  where $R_i(t) \in K[t]$.  Let us
  write $f=\frac pq$, $g = \frac uv$ and $R_n(t) =
  \alpha(t-\lambda_1)\cdots(t-\lambda_m)$.  Then
  $$\sum_{i=0}^{n} R_i\Big(\frac pq\Big)\left(\frac uv\right)^i=0,
  \text{ hence } \sum_{i=0}^{n} R_i\Big(\frac pq\Big) u^iv^{n-i}=0.$$
  By multiplication by $q^d$ for $d = \max \{\deg R_i\}$ (in order
  that $q^d R_i(\frac pq)$ are polynomials) we obtain
  $$q^d R_n\Big(\frac pq\Big) u^n = v \left( -q^dR_{n-1}\Big(\frac
    pq\Big)u^{n-1}-\cdots\right).$$
  As $\gcd(u,v)=1$ then $v$ divides
  the polynomial $q^d R_n(\frac pq)$, then $v$ divides
  $q^{d-m}(p-\lambda_1q)\cdots(p-\lambda_m q)$.  Then all irreducible
  factors of $v$ divide $q$ or $p-\lambda_iq$, $i=1,\ldots,m$.

  Let $\lambda \notin \{\infty,\lambda_1,\ldots,\lambda_m\}$.  Let
  $V_\lambda$ be an irreducible component of $p-\lambda q$, then
  $V_\lambda \cap Z(v)$ is zero dimensional (or empty). Hence $v$ is
  not identically equal to $0$ on $V_\lambda$.  Then for all but
  finitely many $(x,y) \in V_\lambda$ we get:
  $$\sum_{i=0}^{n} R_i(\lambda) g(x,y)^i=0.$$
  Therefore $g$ can only
  reach a finite number of values $c_1,\ldots,c_n$ (the roots of
  $\sum_{i=0}^{n} R_i(\lambda) t^i$). Since $V_\lambda$ is
  irreducible, $g$ is constant on $V_\lambda$.
  
\item (\ref{it:liii}) $\Rightarrow$ (\ref{it:liv}). Clear.
  
\item (\ref{it:liv}) $\Rightarrow$ (\ref{it:li}).  We first give a
  proof that if $g$ is constant along an irreducible component
  $V_\lambda$ of $(p-\lambda q=0)$ then $D_f(g)=0$ on $V_\lambda$ (we
  suppose that $V_\lambda$ is not in the poles of $g$). Let
  $(x_0,y_0)\in V_\lambda$ and $t\mapsto p(t)$ be a local
  parametrization of $V_\lambda$ around $(x_0,y_0)$.  By definition
  of $p(t)$ we have $f(p(t)) = \lambda$, this implies that:
  $$
  \left\langle \frac{dp}{dt}\mid \overline {\grad f} \right\rangle
  = \frac{d(f(p(t))}{dt} = 0$$
  and by hypotheses $g$ is constant on
  $V_\lambda$ this implies $g(p(t))$ is constant and again:
  $$
  \left\langle \frac{dp}{dt}\mid \overline {\grad g} \right\rangle
  = \frac{d(g(p(t))}{dt} = 0.$$
  Then $\grad f$ and $\grad g$ are
  orthogonal around $(x_0,y_0)$ on $V_\lambda$ to the same vector, as
  we are in dimension $2$ this implies that the determinant of
  Jacobian matrix of $(f,g)$ is zero around $(x_0,y_0)$ on
  $V_\lambda$. By extension $D_f(g) =0$ on $V_\lambda$.
  
  We now end the proof: If $g$ is constant on infinitely many
  irreducible components $V_\lambda$ of $(p-\lambda q=0)$ this implies
  that $D_f(g)=0$ on infinitely many $V_\lambda$. Then $D_f(g) =0$ in
  $K(x,y)$.

\end{itemize}
\end{proof}

\subsection{Group of the divisors}

Let $f = \frac pq$, let $\lambda_1,\ldots,\lambda_n \in \hat K$, we
denote by $G(f;\lambda_1,\ldots,\lambda_n)$ the multiplicative group
generated by all the divisors of the polynomials $p-\lambda_i q$,
$i=1,\ldots,n$.

Let
$$d(f) = (\deg f)^2+\deg f.$$

\begin{lemma}
\label{lem:degree}
Let $F_1,\ldots,F_r \in G(f;\lambda_1,\ldots,\lambda_n)$. If $r\ge
d(f)$ then there exists a collection of integers $m_1,\ldots, m_r$
(not all equal to zero) such that
$$g= \prod_{i=1}^r F_i^{m_i} \in C_f.$$
\end{lemma}

\begin{proof}
  Let $\mu \notin \{ \lambda_1,\ldots,\lambda_n \}$, and let $S$ be an
  irreducible component of $(p-\mu q=0)$. Let $\bar S$ be the
  projective closure of $S$.   The
  functions $F_i$ restricted to $\bar S$ have their poles and zeroes
  on the points at infinity of $S$ or on the intersection $S \cap
  Z(F_i) \subset Z(p) \cap Z(q)$.
  
  Let $n : \tilde S \rightarrow \bar S$ be a normalization of $\bar
  S$.  The inverse image under normalisation of the points at infinity
  are denoted by $\{ \gamma_1,\ldots,\gamma_k\}$, their number
  verifies $k \le \deg S \le \deg f$.
  
  At a point $\delta \in Z(p) \cap Z(q)$, the number of points of
  $n^{-1}(\delta)$ is the local number of branches of $S$ at $\delta$
  then it is less or equal than $\ord_\delta (S)$, where $\ord_\delta
  (S)$ denotes the order (or multiplicity) of $S$ at $\delta$ (see
  e.g. \cite{Sh}, paragraph II.5.3).  Then
  \begin{align*}
  \# n^{-1}(\delta) &\le \ord_\delta (S) \le \ord_\delta Z(p- \mu q) \le \ord_\delta Z(p-\mu q) \cdot
  \ord_\delta Z(p) \\ 
    &\le \mult_\delta (p-\mu q,p) = \mult_\delta (p,q) 
  \end{align*}
  where $\mult_\delta
  (p,q)$ is the intersection multiplicity (see e.g. \cite{Fu}).
  Then by B\'ezout theorem:
  $$\sum_{\delta \in Z(p) \cap Z(q)} \# n^{-1}(\delta) \le \sum_{\delta
    \in  Z(p) \cap Z(q)}{\mult_\delta (p,q)} \le \deg p \cdot \deg q \le (\deg f)^2.$$

  Then the inverse image under normalisation of $\cup_{i=1}^r
  S \cap Z(F_i)$ denoted by $\{\gamma_{k+1},\ldots,\gamma_\ell\}$ have
  less or equal than $(\deg f)^2$ elements.  Notice that $\ell \le
  \deg f + (\deg f)^2 = d(f)$.
  
  Now let $\nu_{ij}$ be the order of $F_i$ at $\gamma_j$
  ($i=1,\ldots,r$; $j=1,\ldots,\ell$).  Consider the matrix $M =
  (\nu_{ij})$.  Because the degree of the divisor $(F_i)$ (seen over
  $\tilde S$) is zero we get $\sum_{j=1}^\ell\nu_{ij}=0$, for
  $i=1,\ldots,r$, that means that columns of $M$ are linearly
  dependent.  Then $\mathrm{rk}\, M < \ell \le d(f)$, by hypothesis $r
  \ge d(f)$, then the rows of $M$ are also linearly dependent.  Let
  $m_1(\mu,S),\ldots,m_r(\mu,S)$ such that $\sum_{i=1}^r m_i(\mu,S)
  \nu_{ij} =0$, $j=1,\ldots,\ell$.
  
  Consider the function $g_{\mu,S} = \prod_{i=1}^r
  F_i^{m_i(\lambda,S)}$.  Then this function is regular and does not
  have zeroes or poles at the points $\gamma_j$, because $\sum_{i=1}^r
  m_i(\mu,S) \nu_{ij} =0$.  Then $g_{\mu,S}$ is constant on $S$.
  
  This construction gives a map $(\mu,S) \mapsto
  (m_1(\mu,S),\ldots,m_r(\mu,S))$ from $K$ to $\Zz^r$.  Since $K$ is
  uncountable, there exists infinitely many $(\mu,S)$ with the same
  $(m_1,\ldots,m_r)$.  Then the function $g = \prod_{i=1}^r F_i^{m_i}$
  is constant on infinitely many components of curves of $(p-\mu q=0)$
  and by Lemma \ref{lem:kernel} this implies $g \in C_f$.
\end{proof}

\subsection{Non-composite rational function}

Let $f = \frac pq$. Let $G(f)$ be the multiplicative group generated
by all divisors of the polynomials $p-\lambda q$ for all $\lambda \in
\hat K$. In fact we have
$$G(f) = \bigcup_{(\lambda_1,\ldots,\lambda_n)\in K^n}
G(f;\lambda_1,\ldots,\lambda_n).$$

\begin{definition}
  A family $F_1,\ldots,F_r \in G(f)$ is \defi{$f$-free} if
  $(m_1,\ldots,m_r)\in \Zz^r$ is such that $ \prod_{i=1}^r F_i^{m_i}
  \in C_f$ then $(m_1,\ldots,m_r) = (0,\ldots,0)$.
  
  A $f$-free family $F_1,\ldots,F_r \in G(f)$ is \defi{$f$-maximal} if
  for all $F\in G(f)$, $\{F_1,\ldots,F_r,F\}$ is not $f$-free.
\end{definition}

\begin{theorem}
\label{th:noncomposite}
Let $f\in K(x,y)$, $\deg f >0$. Then the following conditions are
equivalent:
\begin{enumerate}
\item \label{it:ti} $\deg f = \min \left\lbrace \deg g \mid g\in C_f
    \setminus K \right\rbrace$;
\item \label{it:tii} $\sigma(f)$ is finite;
\item \label{it:tiii} $C_f = K(f)$;
\item \label{it:tiv} $f$ is non-composite.
\end{enumerate}
\end{theorem}

\begin{remark}
  This does not give a new proof of ``$\sigma(f)$ is finite $\Leftrightarrow f$ is non-composite'' because we use Bertini-Krull
  theorem.
\end{remark}

\begin{remark}
  The proof (\ref{it:ti}) $\Rightarrow$ (\ref{it:tii}) is somewhat
  easier than in \cite{St}, whereas (\ref{it:tii}) $\Rightarrow$
  (\ref{it:tiii}) is more difficult.
\end{remark}

\begin{proof}
  \ 
\begin{itemize}
\item (\ref{it:ti}) $\Rightarrow$ (\ref{it:tii}).  Let us suppose that
  $\sigma(f)$ is infinite. Set $f = \frac pq$, with $\gcd(p,q)=1$. For all
  $\alpha \in \sigma(f)$, let $F_\alpha$ be an irreducible divisor of
  $p-\alpha q$, such that $\deg F_\alpha < \deg f$.  By
  Lemma \ref{lem:degree} there exists a $f$-maximal family $\{
  F_1,\ldots,F_r\}$   with $r\le d(f)$.  Moreover $r\ge 1$
  because $\{ F_\alpha \}$ is $f$-free: if not there exists $k\neq 0$
  such that $F_\alpha^k \in C_f$ then $F_\alpha\in C_f$, but $\deg
  F_\alpha < \deg f$ that contradicts the hypothesis of minimality.
  
  Now the collection $\{ F_1,\ldots,F_r,F_\alpha\}$ is not $f$-free,
  so that there exist integers $\{
  m_1(\alpha),\ldots,m_r(\alpha),m(\alpha)\}$, with $m(\alpha)\neq 0$,
  such that
  $$F_1^{m_1(\alpha)}\cdots F_r^{m_r(\alpha)}\cdot
  F_\alpha^{m(\alpha)} \in C_f.$$
  Since $\sigma(f)$ is infinite then is equal to $\hat K$ minus a finite number of values
(see Theorem \ref{th:composite}) then $\sigma(f)$ is uncountable and the map
  $\alpha \mapsto ( m_1(\alpha),\ldots,m_r(\alpha),m(\alpha))$ is not
  injective.  Let $\alpha \neq \beta$ such that
  $m_i(\alpha)=m_i(\beta) = m_i$, $i=1,\ldots,r$ and
  $m(\alpha)=m(\beta)=m$.  Then $F_1^{m_1}\cdots F_r^{m_r}\cdot
  F_\alpha^{m} \in C_f$ and $F_1^{m_1}\cdots F_r^{m_r}\cdot
  F_\beta^{m} \in C_f$, it implies that $(F_\alpha/F_\beta)^m \in
  C_f$, therefore $F_\alpha/F_\beta \in C_f$.
  
  Now $\deg \frac{F_\alpha}{F_\beta} < \deg f$, then by the hypothesis
  of minimality it proves $\frac{F_\alpha}{F_\beta}$ is a constant.
  Let $a\in K^*$ such that $F_\alpha = a F_\beta$, by definition
  $F_\alpha$ divides $p-\alpha q$, but moreover $F_\alpha$ divides
  $p-\beta q$ (as $F_\beta$ do).  Then as $F_\alpha$ divides both
  $p-\alpha q$ and $p-\beta q$, $F_\alpha$ divides $p$ and $q$, that
  contradicts $\gcd(p,q)=1$.

\item (\ref{it:tii}) $\Rightarrow$ (\ref{it:tiii}).  Let $f=\frac pq$,
  $\sigma(f)$ finite and $g \in C_f$, we aim at proving that $g\in
  K(f)$.  The proof will be done in several steps:
\begin{enumerate}
\item[(a)] \emph{Reduction to the case $g=\frac{u}{q^\ell}$.}  Let $g
  = \frac uv \in C_f$, then $f$ and $g$ are algebraically dependent,
  then there exists a polynomial in $f$ and $g$ that vanishes. As
  before let us write
  $$\sum_{i=0}^{n} R_i(f)g^i=0$$
  where $R_i(t) \in K[t]$.  As $f =
  \frac pq$, $g=\frac uv$ then
  $$\sum_{i=0}^{n} R_i\Big(\frac pq\Big)\left(\frac uv\right)^i=0,
  \text{ hence } \sum_{i=0}^{n} R_i\Big(\frac pq\Big) u^iv^{n-i}=0.$$
  By multiplication by $q^d$ for $d = \max \{\deg R_i\}$ (in order
  that all $q^d R_i(\frac pq)$ are polynomials) we get:
  $$q^d R_n\Big(\frac pq\Big) u^n = v \left( -q^dR_{n-1}\Big(\frac
    pq\Big)u^{n-1}-\cdots\right).$$
  As $\gcd(u,v)=1$ then $v$ divides
  the polynomial $q^d R_n(\frac pq)$; we write $vu' = q^d R_n(\frac
  pq)$ then
  $$g = \frac uv = \frac{uu'}{q^d R_n(\frac pq)}.$$
  But $R_n(\frac pq)
  \in K(\frac pq)$ then $\frac{uu'}{q^d} \in C_f$, but also we have
  that $g\in K(f)$ if and only if $\frac{uu'}{q^d} \in K(f)$.  This
  proves the reduction.

\item[(b)] \emph{Reduction to the case $g=qu$.}  Let $g = \frac
  {u}{q^\ell} \in C_f$, $\ell \ge 0$.  As 
  $\sigma(f)$ is finite by Lemma
  \ref{lem:kernel} we choose $\lambda \in K$ such that $p-\lambda q$
  is irreducible and $g \in C_f$ is constant (equal to $c$) on
  $p-\lambda q$.  As $g=\frac {u}{q^\ell}$,we have $p- \lambda q$
  divides $u-cq^\ell$. We can write:
  $$u-c q^\ell = u'(p-\lambda q).$$
  Then
  $$\frac{u}{q^\ell} = \frac{u'}{q^{\ell-1}} \Big( \frac pq -
  \lambda\Big) + c.$$
  As $\frac{u}{q^\ell}$ and $f=\frac pq$ are in
  $C_f$ we get $ \frac{u'}{q^{\ell-1}}\in C_f$; moreover
  $\frac{u}{q^\ell} \in K(f)$ if and only if $
  \frac{u'}{q^{\ell-1}}\in K(f)$. By induction on $\ell \ge 0$ this
  prove the reduction.
  
\item[(c)] \emph{Reduction to the case $g=q$.}  Let $g=qu \in C_f$.
  $g$ is constant along the irreducible curve $(p-\lambda q=0)$.  Then
  $qu=u_1(p-\lambda q)+c_1$.
  
  Let $\deg p = \deg q$. Then $q^hu^h=u_1^h(p^h-\lambda q^h)$ (where
  $P^h$ denotes the homogeneous part of higher degree of the
  polynomial $P$). Then $p^h-\lambda q^h$ divides $q^hu^h$ for
  infinitely many $\lambda \in K$.  As $\gcd(p,q)=1$ this gives a
  contradiction.
  
  Hence $\deg p \neq \deg q$. We may assume $\deg p > \deg q$
  (otherwise $qu\in C_f$ and $\frac pq \in C_f$ implies $pu \in C_f$).
  Then we write:
  $$qu=qu_1\Big( \frac pq - \lambda \Big)+c_1,$$
  that proves that
  $qu_1 \in C_f$ and that $qu \in K(f)$ if and only if $qu_1 \in
  K(f)$.  The inequality $\deg p > \deg q$ implies that $\deg u_1 <
  \deg u$.  We continue by induction, $qu_1= qu_2(\frac pq
  -\lambda)+c_2$, with $\deg u_2 < \deg u_1$,..., until we get $\deg
  u_n = 0$ that is $u_n\in K^*$.  Thus we have prove firstly that
  $qu_n \in C_f$, that is to say $q\in C_f$, and secondly that $qu \in
  K(f)$ if and only if $q\in K(f)$.

\item [(d)] \emph{Case $g=q$.}  If $q\in C_f$ then $q$ is constant
  along the irreducible curve $(p-\lambda q=0)$ then $q= a(p-\lambda
  q)+c$, $a\in K^*$. Then
  $$q=\frac{c}{1-a(\frac pq - \lambda)} \in K\Big(\frac pq\Big) =
  K(f).$$

\end{enumerate}

\item (\ref{it:tiii}) $\Rightarrow$ (\ref{it:tiv}).  Let us assume
  that $C_f= K(f)$ and that $f$ is composite, then there exist $r\in
  K(t)$, $\deg r \ge 2$ and $g \in K(x,y)$ such that $f=r\circ g$. By
  the formula $\deg f = \deg r \cdot \deg g$ we get $\deg f > \deg g$.
  Now if $r = \frac ab$ then we have a relation $b(g)f=a(g)$, then $f$
  and $g$ are algebraically dependent, hence by Lemma
  \ref{lem:kernel}, $g\in C_f$. As $C_f=K(f)$, there exists $s\in
  K(t)$ such that $g = s \circ f$. Then $\deg g \ge \deg f$.  That
  yields to a contradiction.

\item (\ref{it:tiv}) $\Rightarrow$ (\ref{it:ti}).  Assume that $f$ is
  non-composite and let $g\in C_f$ of minimal degree.  By Corollary
  \ref{cor:kernel} we get $C_f=C_g$, then $\deg g = \min \left\lbrace
    \deg h \mid h\in C_g \setminus K \right\rbrace$.  Then by the
  already proved implication (\ref{it:ti}) $\Rightarrow$
  (\ref{it:tiii}) for $g$, we get $C_g=K(g)$. Then $f\in
  C_f=C_g=K(g)$, then there exists $r\in K(t)$ such that $f = r\circ
  g$, but as $f$ is non-composite then $\deg r=1$, hence $\deg f=\deg
  g = \min \left\lbrace \deg h \mid h\in C_f \setminus K
  \right\rbrace$.

\end{itemize}
\end{proof}

\section{Order of reducibility of rational functions in two variables}
\label{sec:order}

Let $f = \frac pq \in K(x,y)$; for all $\lambda \in \hat K$, let
$n_\lambda$ be the number of irreducible components of $p-\lambda q$.
Let
$$\rho(f) = \sum_{\lambda \in \hat K} (n_\lambda -1).$$

By Theorem \ref{th:composite}, $\rho(f)$ is finite if and only if $f$
is non-composite. We give a bound for $\rho(f)$.  Recall that we
defined:
$$d(f) = (\deg f)^2+\deg f.$$

\begin{theorem}
\label{th:twovar}
Let $K$ be an algebraic closed field of characteristic $0$.  
If $f \in K(x,y)$ is non-composite then
$$\rho(f) < d(f).$$
\end{theorem}

\begin{proof}
  First notice that $K$ can be supposed uncountable, otherwise it can be 
  embedded into an uncountable field $L$ and the spectrum in $K$ would 
  be included in the spectrum in $L$.
  
  Let us assume that $f$ is non-composite, then by Theorem
  \ref{th:composite} and its corollary we have that $\sigma(f)$ is
  finite: $\sigma(f)=\{\lambda_1,\ldots,\lambda_r\}$.  We suppose that
  $\rho(f) \ge d(f)$.  Let $f = \frac pq$. We decompose the
  polynomials $p-\lambda_i q$ in irreducible factors, for
  $i=1,\ldots,r$:
  $$p-\lambda_i q = \prod_{j=1}^{n_i} F_{i,j}^{k_{i,j}},$$
  where $n_i$
  stands for $n_{\lambda_i}$.  Notice that since $\gcd(p,q) =1$ then
  $F_{i,j}$ divides $p-\lambda_i q$ but do not divides any of $p-\mu
  q$, $\mu \neq \lambda_i$.  The collection $\left\lbrace
    F_{1,1},\ldots,F_{1,n_1-1},\ldots ,
    F_{r,1},\ldots,F_{r,n_r-1}\right\rbrace$, is included in
  $G(f,\lambda_1,\ldots,\lambda_r)$ and contains $\rho(f) \ge d(f)$
  elements, then Lemma \ref{lem:degree} provides a collections $\{
  m_{1,1},\ldots,m_{1,n_1-1},\ldots , m_{r,1},\ldots,m_{r,n_r-1}\}$ of
  integers (not all equal to $0$) such that
\begin{equation}
g = \prod_{i=1}^r \prod_{j=1}^{n_i-1} F_{i,j}^{m_{i,j}} \in C_f.
\label{eq:defg}
\end{equation}

By Theorem \ref{th:noncomposite} it implies that $g \in K(f)$, then
$g=\frac{u(f)}{v(f)}$, where $u,v \in K[t]$.  Let $\mu_1,\ldots,\mu_k$
be the roots of $u$ and $\mu_{k+1},\ldots,\mu_\ell$ the roots of $v$.
Then
$$g =\frac{u(\frac pq)}{v(\frac pq)} = \alpha \frac{\prod_{i=1}^k
  \frac pq - \mu_i}{\prod_{i=k+1}^\ell \frac pq - \mu_i}$$
so that
\begin{equation}
g = \alpha q^{\ell-2k}\frac{\prod_{i=1}^k  p - \mu_i q}{\prod_{i=k+1}^\ell  p - \mu_i q}.
\label{eq:fracg}
\end{equation}

If $m_{i_0,j_0} \neq 0$ then by the definition of $g$ by equation
(\ref{eq:defg}) and by equation (\ref{eq:fracg}), we get that
$F_{i_0,j_0}$ divides one of the $p-\mu_i q$ or divides $q$.  If
$F_{i_0,j_0}$ divides $p - \mu_i q$ then $\mu_i= \lambda_{i_0} \in
\sigma(f)$.  If $F_{i_0,j_0}$ divides $q$ then $\mu_i=\infty$, so that
$\infty \in \sigma(f)$.  In both cases $p - \lambda_{i_0} q$ appears
in formula (\ref{eq:fracg}) at the numerator or at the denominator of
$g$.  Then $F_{i_0,n_{i_0}}$ should appears in decomposition
(\ref{eq:defg}), that gives a contradiction.  Then $\rho(f) < d(f)$.
\end{proof}

\section{Extension to several variables}
\label{sec:several}

We follows the lines of the proof of \cite{Na2}.  We will need a
result that claims that the irreducibility and the degree of a family
of polynomials remain constant after a generic linear change of
coordinates.  For $\xubar = (x_1,\ldots,x_n)$ and a matrix $B =
(b_{ij}) \in Gl_n(K)$, we denote the new coordinates by $B\cdot
\xubar$:
$$B\cdot \xubar = (\sum_{j=1}^n b_{1j} x_j,\ldots, \sum_{j=1}^n b_{nj}
x_j).$$

\begin{proposition}
\label{prop:changecoord}
Let $K$ be an infinite field.  Let $n\ge 3$ and $p_1,\ldots,p_\ell \in
K[x_1,\ldots,x_n]$ be irreducible polynomials. Then there exists a
matrix $B \in Gl_n(K)$ such that for all $i=1,\ldots,\ell$ we get:
\begin{itemize}
\item $p_i(B\cdot \xubar)$ is irreducible in
  $\overline{K(x_1)}[x_2,\ldots,x_n]$;
\item $\deg_{(x_2,\ldots,x_n)} p_i(B\cdot \xubar) =
  \deg_{(x_1,\ldots,x_n)} p_i$.
\end{itemize}
\end{proposition}

The proof of this proposition can be derived from \cite[Ch. 5, Th.
3D]{Scm} or by using \cite[Prop. 9.31]{FJ}. See \cite{Na2} for
details.

\bigskip

Now we return to our main result.
\begin{theorem}
  Let $K$ be an algebraically closed field of characteristic $0$.
  Let $f \in K(\xubar)$ be non-composite then $\rho(f) <
  (\deg f)^2+\deg f$.
\end{theorem}

\begin{proof}
  We will prove this theorem by induction on the number $n$ of
  variables.  For $n=2$, we proved in Theorem \ref{th:twovar} that
  $\rho(f) < (\deg f)^2+\deg f$.
  
  Let $f= \frac pq \in K(\xubar)$, with $\xubar = (x_1,\ldots,x_n)$.
  We suppose that $f$ is non-composite.  For each $\lambda \in
  \sigma(f)$ we decompose $p-\lambda q$ into irreducible factors:
\begin{equation}
p-\lambda q = \prod_{i=1}^{n_\lambda} F_{\lambda,i}^{r_{\lambda,i}}.
\label{eq:factor}
\end{equation}

We fix $\mu \notin \sigma(f)$.  We apply Proposition
\ref{prop:changecoord} to the polynomials $p-\mu q$ and
$F_{\lambda,i}$, for all $\lambda \in \sigma(f)$ and all
$i=1,\ldots,n_\lambda$.  Then the polynomials $p(B\cdot\xubar)-\mu
q(B\cdot\xubar)$ and $F_{\lambda,i}(B\cdot\xubar)$ are irreducible in
$\overline{K(x_1)}[x_2,\ldots,x_n]$ and their degrees in
$(x_2,\ldots,x_n)$ are equals to the degrees in $(x_1,\ldots,x_n)$ of
$p-\mu q$ and $F_{\lambda,i}$.

Let denote by $k = \overline {K(x_1)}$.  This is an uncountable field,
algebraically closed of characteristic zero.  Now $p(B\cdot\xubar)-\mu
q(B\cdot\xubar)$ is irreducible, then $f(B\cdot\xubar)$ is
non-composite in $k(x_2,\ldots,x_n)$.

Now equation (\ref{eq:factor}) become:
$$p(B\cdot\xubar)-\lambda q(B\cdot\xubar) = \prod_{i=1}^{n_\lambda}
F_{\lambda,i}(B\cdot\xubar)^{r_{\lambda,i}}.$$
Which is the
decomposition of $p(B\cdot\xubar)-\lambda q(B\cdot\xubar)$ into
irreducible factors in $k(x_2,\ldots,x_n)$.  Then
$$\sigma(f) \subset \sigma(f(B\cdot\xubar)),$$
where $\sigma(f)$ is a
subset of $K$, and $\sigma(f(B\cdot\xubar))$ is a subset of
$k=\overline{K(x_1)}$.  As $n_\lambda$ is also the number of distinct
irreducible factors of $p(B\cdot\xubar)-\lambda q(B\cdot\xubar)$ we
get:
$$\rho(f) \le \rho(f(B\cdot\xubar)).$$

Now suppose that the result is true for $n-1$ variables.  Then for $f
(B\cdot\xubar) \in k(x_2,\ldots,x_n)$ we get:
$$
\rho(f(B\cdot\xubar)) < (\deg_{(x_2,\ldots,x_n)}
f(B\cdot\xubar))^2+(\deg_{(x_2,\ldots,x_n)} f(B\cdot\xubar)).$$
Hence:
\begin{align*}
  \rho(f) &\le  \rho(f(B\cdot\xubar)) \\
  &<  (\deg_{(x_2,\ldots,x_n)} f(B\cdot\xubar))^2+(\deg_{(x_2,\ldots,x_n)} f(B\cdot\xubar)) \\
  &=  (\deg_{(x_1,\ldots,x_n)}f)^2+(\deg_{(x_1,\ldots,x_n)}f) \\
  &=  (\deg f)^2+(\deg f) \\
\end{align*}
\end{proof}

If for $n=2$ we start the induction with Lorenzini's bound $\rho(f) <
(\deg f)^2$ we obtain with the same proof the following result for
several variables, for $K$ of any characteristic $K$ and a better
bound:
\begin{theorem}
  Let $K$ be an algebraically closed field.  Let $f \in K(\xubar)$ be
  non-composite then $\rho(f) < (\deg f)^2$.
\end{theorem}


\end{document}